\documentclass[10]{amsart}
\usepackage[T1]{fontenc}
\usepackage{latexsym}
\usepackage{amssymb}
\usepackage{amsfonts}
\usepackage{float}
\usepackage[all]{xy}
\theoremstyle{remark}

\theoremstyle{definition}

\numberwithin{equation}{section}

\newcommand{\R}{\text{\bf R}}           

\newcommand{\Cal}{\mathcal}

\newcommand{\be}{\begin{equation}}

\newcommand{\beu}{\begin{equation*}}

\newcommand{\IP}[2]{\langle#1\,, #2\rangle}     

\begin{document}

\baselineskip=19.5pt

\sloppy

\title[Cohomologie $L^{p}$]
{Cohomologie $L^{p}$ et formes harmoniques}

\author{No\"el Lohou\'e}

\thanks{{\em Mots cl\'es :} cohomology $L^{2}$, d\'ecomposition de de Rham, forme harmonique $L^{p}$.}
\address{CNRS and D\'epartement de math\'ematiques, Universit\'e Paris 11, B\^atiment 425, 91405 Orsay, France}
\email{Noel.Lohoue@math.u-psud.fr} \maketitle
\begin{abstract}
{On montre que si une vari\'et\'e riemannienne admet un rev\^etement universel \`a g\'eom\'etrie born\'ee et si $0$ est hors du spectre ou un point isol\'e du spectre du laplacien 
sur les formes de degr\'e $\ell$ alors il existe $1<p<2$ tel que pour tout $p<r<p^{\prime}$ la d\'ecomposition de Hodge - de Rham pour les formes $L^{r}$ soit vraie 
($p^{\prime}$ d\'esigne le conjugu\'e de $p$).}
\end{abstract}

\section{Introduction.}

Soit $M$ une variété riemannienne, complète de dimension $n$, avec un élément de volume $d\sigma$.

Si $\ell$ est un entier, $0\leq\ell\leq n$, on note $C_{0}^{\infty}(\Lambda^{\ell}M)$ l'espace des $\ell$-formes diff\'erentielles sur $M$, de classe $C^{\infty}$, \`a support compact. 
Pour tous $\omega\in C_{0}^{\infty}(\Lambda^{\ell}M)$ et $x\in M$ on d\'esigne par $\mid\mid\omega_{x}\mid\mid$ la norme de $\omega_{x}$ induite par la structure riemannienne. 

Si $1\leq p<\infty$, soit $L^{p}(\Lambda^{\ell}M)$ le compl\'et\'e de $C_{0}^{\infty}(\Lambda^{\ell}M)$ pour la norme
$$\mid\mid\omega\mid\mid_{p}^{p}=\int_{M}\mid\mid\omega_{x}\mid\mid^{p}d\sigma(x).$$
Si $\omega$ est une forme mesurable et presque partout born\'ee on pose
$$\mid\mid\omega\mid\mid_{\infty}=esssup_{x\in M}\mid\mid\omega_{x}\mid\mid.$$
Pour $p>1$ on pose $p^{\prime}=p/(p-1)$. Si $T$ est une application lin\'eaire de $L^{p}(\Lambda^{\ell}M)$ dans $L^{q}(\Lambda^{\ell}M)$, on note 
$\mid\mid T\mid\mid_{p\rightarrow q}$ la norme de l'op\'erateur $T$.

Rappelons que $M$ est {\it \`a g\'eom\'etrie born\'ee} si pour tout $x\in M$ il existe un diff\'eomorphisme $\psi_{x}$ de la boule unit\'e de $\R^{n}$ sur la boule unit\'e de $M$ centr\'ee en $x$ 
avec un contr\^ole uniforne des d\'eriv\'ees de $\psi_{x}$. 

D\'esignons par $d$ l'op\'erateur de d\'erivation ext\'erieure sur $M$ et $\delta$ son adjoint au sens du produit scalaire sur $L^{2}(\Lambda^{\ell}M)$ induit par la mesure $d\sigma$.

Soit $\delta=d\delta+\delta d$ le laplacien de Hodge - de Rham sur $M$. D'apr\`es un th\'eor\`eme de de Rham toute forme $\omega$ de carr\'e sommable ($\omega\in L^{2}(\Lambda^{\ell}M)$), 
de degr\'e $\ell$, s'\'ecrit de fa\c con unique comme
$$\omega=\omega_{1}+\omega_{2}+\omega_{3}$$
o\`u $\omega_{1}$ est adh\'erent \`a l'image $d C_{0}^{\infty}(\Lambda^{\ell-1}M)$, $\omega_{2}$ est adh\'erent \`a l'image $\delta C_{0}^{\infty}(\Lambda^{\ell+1}M)$, et o\`u $\omega_{3}$ est harmonique ($\Delta\omega_{3}=0$). Par la suite $c$ d\'esignera une constante positive dont la valeur nous importe peu et $\tau$ est une constante qui intervient pour la premi\`ere fois \`a la page $3$ et sera la m\^eme par la suite.

On d\'esignera par ${\Cal H}_{p}$ l'espace des formes harmoniques de $L^{p}(\Lambda^{\ell}(M)$.

Le but de cette note est de g\'en\'eraliser, dans la mesure du possible, cette d\'ecomposition \`a $L^{p}(\Lambda^{\ell}(M)$, pour $p\neq 2$. Plus pr\'ecis\'ement on prouvera l'\'enonc\'e suivant :

{\bf Théorème :} {\em Soit $M$ une variété riemmenienne complète. On suppose que
\begin{itemize}
\item[(i)] $M$ admet une revêtement $\tilde{M}$ à géométrie bornée.
\item[(ii)] $0$ est un point isolé du spectre de $\Delta$ sur $L^{2}(\Lambda^{\ell}M)$ ou $0$ n'appartient pas au spectre de $\Delta$.
\end{itemize}
Alors il existe des r\'eels $p_{1}$ et $p_{2}$, $1<p_{1}<2<p_{2}$ tels que pour tout $p_{1}<p<p_{2}$, toute forme $\omega\in L^{p}(\Lambda^{\ell}M)$
s'écrive de façon unique : $$\omega=d\omega_{1}+\delta\omega_{2}+\omega_{3}$$ 
avec}
$$\Delta\omega_{3}=0,\;\;\mid\mid\omega_{i}\mid\mid_{p}\leq c\mid\mid\omega\mid\mid_{p}\;\;\; (i=1,2,3)$$
$$\mid\mid d\omega_{1}\mid\mid_{p}\leq c\mid\mid\omega\mid\mid_{p},\;\;\;\mid\mid\delta\omega_{2}\mid\mid_{p}\leq c\mid\mid\omega\mid\mid_{p}.$$ 

{\bf Remarque.} Un rapporteur nous a signal\'e que Y. A. Kordyukov a obtenu des r\'esultats dans cette direction. Il a aussi \'etabli des estimations \`a priori sur les vari\'et\'es riemanniennes. 
Il convient de signaler que de telles estimations \`a priori ont \'et\'e \'etablies dans \cite{Lo06}. Elles dates d'ailleurs de 1984 (preprint d'Orsay). Elles ne sont pas utilisables ici car $M$ n'est pas 
suppos\'ee \`a g\'eom\'etrie born\'ee.

Le th\'eor\`eme \'enonc\'e ci-dessus a pour origine une question de P. Pansu. Les deux rapporteurs nous ont signal\'e des coquilles dans le texte inital; nous les en remer\c cions.

\section{Démonstration du théorème.}

D'après la décomposition de de Rham-Hodge
\beu
L^{2}(\Lambda^{\ell}M)=\overline{dC^{\infty}_{0}(\Lambda^{\ell-1}M)}\oplus\overline{\delta C^{\infty}_{0}(\Lambda^{\ell+1}M)}\oplus{\Cal H}_{2}
\end{equation*}
où $C^{\infty}_{0}(\Lambda^{\ell}M)$ est l'espace des formes différentielles $C^{\infty}$ à support compact, et o\`u ${\Cal H}_{2}$ est l'espace des formes
harmonique de carré intégrable. On note $H$ la projection orthogonale sur ${\Cal H}_{2}$ alors il est clair que si $P_{t}^{\ell}$ désigne la
solution fondamentale de l'équation de la chaleur sur les $\ell$-formes, alors pour toute forme $\omega\in L^{2}(\Lambda^{\ell}M)$, on a :
\beu
H\omega=\lim_{t\rightarrow +\infty}P_{t}^{\ell}\omega.
\end{equation*}

{\bf Id\'ee de la preuve du th\'eor\`eme.} On commence par montrer que si $H$ est le projecteur sur la composante harmonique de $L^{2}(\Lambda^{\ell}(M))$, il existe $1<q_{0}<2$ tel que 
pour tout $q_{0}<p<q_{0}^{\prime}$ il existe une constante $c$ telle que 
$$\mid\mid H\omega\mid\mid_{p}\leq c\mid\mid\omega\mid\mid_{p}.$$
Dans une seconde \'etape, on montre que :
$$\mid\mid\int_{0}^{\infty}P_{t}^{\ell}(1-H)\omega dt\mid\mid_{p}\leq c\mid\mid\omega\mid\mid_{p},\;\;\;p_{1}<p<p_{2}$$
o\`u $P_{t}^{\ell}$ est le noyau de la chaleur sur les formes $\ell$-formes. Soit $\omega_{0}=(1-H)\omega$. On veut voir que si l'on \'ecrit 

$\omega=\Delta G\omega+H\omega=\delta d G(1-H)\omega+d\delta G(1-H)\omega+H\omega$, $G(1-H)=\int_{0}^{\infty}P_{t}^{\ell}(1-H)\omega dt$ alors
$$\mid\mid\delta dG(1-H)\omega\mid\mid_{p}\leq c\mid\mid\omega\mid\mid_{p}\;\;\text{ et }\;\;\mid\mid d\delta G(1-H)\omega\mid\mid_{p}\leq c\mid\mid\omega\mid\mid_{p}.$$
Par exemple, si $\omega_{0}=(1-H)\omega$
\begin{eqnarray*}
d\delta G\omega_{0}&=&d\delta G(1-H)\omega\\
&=&d\delta\Delta^{-\frac{1}{2}}\Delta^{-\frac{1}{2}}\omega_{0}\\
&=&\Delta^{-\frac{1}{2}}d\delta\Delta^{-\frac{1}{2}}\omega_{0}\\
&=&\Delta^{-\frac{1}{2}}d\delta\Delta^{-\frac{1}{2}}\omega_{0}
\end{eqnarray*}
et
$$\mid\mid d\delta G\omega_{0}\mid\mid_{p}\leq\mid\mid\Delta^{-\frac{1}{2}}d\mid\mid_{p\rightarrow p}\mid\mid\delta\Delta^{-\frac{1}{2}}\omega_{0}\mid\mid_{p}\leq c\mid\mid\omega\mid\mid_{p}.$$

Pour prouver le th\'eor\`eme, on aura besoin du lemme suivant :

{\bf Lemme 1.} {\em Sous les hypoth\`eses du th\'eor\`eme, il existe $q_{0}<2$ tel que, pour tout $q_{0}<p<q_{0}^{\prime}$}, 
\beu
\mid\mid H\omega\mid\mid_{p}\leq c\mid\mid\omega\mid\mid_{p}.
\end{equation*}

{\bf Preuve du lemme :} L'hypoth\`ese (ii) montre qu'il existe $c>0$ v\'erifiant $\IP{\Delta\omega}{\omega}\geq c\mid\mid\omega\mid\mid_{2}^{2}$ pour toute forme 
orthogonale \`a ${\Cal H}_{2}$. On note ${\Cal H}_{2}^{\perp}$ le suppl\'ementaire orthogonal de ${\Cal H}_{2}$; alors la norme de $P_{t}^{\ell}$ restreinte à ${\Cal H}_{2}^{\perp}$ satisfait :
\beu
\mid\mid P_{t}^{\ell}\mid_{{\Cal H}^{\perp}}\mid\mid_{2\rightarrow 2}\leq ce^{-\tau_{0} t}\;\text{ avec }\tau_{0}>0.
\end{equation*}
Si l'on écrit $\omega=\omega^{\prime}+\omega^{\prime\prime}$, avec $\omega^{\prime}\in{\Cal H}_{2}^{\perp}$ et $\omega^{\prime\prime}\in{\Cal H}_{2}$, on voit que :
\begin{itemize}
\item[] $\mid\mid P_{t}^{\ell}\omega^{\prime}\mid\mid_{2}\leq ce^{-\tau_{0} t}\mid\mid\omega^{\prime}\mid\mid_{2}\rightarrow 0$ quand $t\rightarrow +\infty$
\item[] $P_{t}^{\ell}\omega^{\prime\prime}=\omega^{\prime\prime}$ et $H\omega=\omega^{\prime\prime}=\lim_{t\rightarrow +\infty}P_{t}^{\ell}\omega$.
\end{itemize}
On veut examiner :
\beu
\omega-H\omega=-\int_{0}^{\infty}\frac{\partial}{\partial s}P_{s}^{\ell}\omega ds.
\end{equation*}
Cette expression a un sens si $\omega$ est $C^{\infty}$ à support compact, car
$\frac{\partial}{\partial s}P_{s}^{\ell}\omega=\frac{\partial}{\partial s}P_{s}^{\ell}\omega^{\prime}$ dans l'écriture précédente; dans la décomposition
spectrale, $\frac{\partial}{\partial s}P_{s}^{\ell}\omega$ correspond à la fonction $\lambda e^{-t\lambda}$ pour $\lambda >c$ et, si $t>t_{0}$, cette
fonction est bornée par $ce^{-\tau_{1}t}$, par conséquent :
\begin{eqnarray*}
\mid\mid -\int_{0}^{\infty}\frac{\partial}{\partial s}P_{s}^{\ell}\omega ds\mid\mid_{2}&\leq&\mid\mid P_{t_{0}}^{\ell}\omega-\omega\mid\mid_{2}+
\mid\mid\int_{t_{0}}^{\infty}\frac{\partial}{\partial s}P_{s}^{\ell}\omega ds\mid\mid_{2}\\
&\leq&\mid\mid\omega\mid\mid_{2}+c\int_{t_{0}}^{\infty}e^{-\tau_{1} t}dt
\mid\mid\omega\mid\mid_{2}.
\end{eqnarray*}
On choisira $t_{0}$ par la suite; nous avons
\beu
\omega-H\omega=-\int_{0}^{t_{0}}\frac{\partial}{\partial s}P_{s}^{\ell}\omega ds-\int_{t_{0}}^{\infty}\frac{\partial}{\partial s}P_{s}^{\ell}\omega ds.
\end{equation*}
Débarrassons nous du terme le moins gênant : pour tout $1\leq p<\infty$,
\beu
\mid\mid\int_{0}^{t_{0}}\frac{\partial}{\partial s}P_{s}^{\ell}\omega ds\mid\mid_{p}=\mid\mid P_{t_{0}}^{\ell}\omega-\omega\mid\mid_{p}
\leq c\mid\mid\omega\mid\mid_{p}.
\end{equation*}
Nous avons utilisé implicitement l'inégalité $\mid\mid P_{t_{0}}^{\ell}\omega\mid\mid_{p}\leq e^{ct_{0}}\mid\mid\omega\mid\mid_{p}$ qu'on peut trouver
dans \cite{Lo06}. Par ailleurs, en utilisant la propri\'et\'e du semigroupe :
\beu
\int_{t_{0}}^{\infty}\frac{\partial}{\partial s}P_{s}^{\ell}\omega ds=
\int_{\frac{t_{0}}{2}}^{\infty}\frac{\partial}{\partial s}P_{s+\frac{t_{0}}{2}}^{\ell}\omega ds=
\int_{\frac{t_{0}}{4}}^{\infty}\Delta P_{\frac{t_{0}}{4}}^{\ell}\circ P_{s+\frac{t_{0}}{4}}^{\ell}\omega ds.
\end{equation*}
On a déjà vu que
\begin{equation}\label{eq1}
\mid\mid\frac{\partial}{\partial s}P_{s+\frac{t_{0}}{2}}^{\ell}\omega\mid\mid_{2}\leq ce^{-\tau_{1} s}\mid\mid\omega\mid\mid_{2};\;\text{ on prend }\tau=\inf\{\tau_{0},\tau_{1}\}
\end{equation}
Montrons que :
\begin{eqnarray}\label{eq2}
\mid\mid\frac{\partial}{\partial s}P_{s+\frac{t_{0}}{4}}^{\ell}\omega\mid\mid_{1}&\leq&
ce^{\alpha t}\mid\mid\omega\mid\mid_{1}\text{ avec }\alpha>0.
\end{eqnarray}
On a :
\begin{eqnarray}
P_{s+\frac{t_{0}}{4}}^{\ell}&=&P_{\frac{t_{0}}{4}}^{\ell}\circ P_{s}^{\ell}\nonumber\\
\frac{\partial}{\partial s}P_{s+\frac{t_{0}}{2}}^{\ell}&=&\Delta P_{\frac{t_{0}}{4}}^{\ell}\circ P_{s+\frac{t_{0}}{4}}^{\ell}\nonumber
\end{eqnarray}
De plus, d'apr\`es \cite{Lo06}, on a l'in\'egalit\'e :
\begin{equation}\label{eq3}
\mid\mid P_{s}^{\ell}\omega\mid\mid_{1}\leq c_{1}e^{\alpha t}\mid\mid\omega\mid\mid_{1}.
\end{equation}
On veut examiner $\Delta P_{\frac{t_{0}}{4}}^{\ell}$ sur $L^{1}(\Lambda^{\ell}M)$. Il suffit de l'examiner sut $L^{\infty}$. Soit $\tilde{\Delta}$
le laplacien sur les $\ell$-formes sur $\tilde{M}$ et $\tilde{P}_{t}^{\ell}$ la solution fondamentale de l'équation de la chaleur associée. One note
$\pi:\tilde{M}\rightarrow M$ le revêtement, de sorte que
\beu
\mid\mid\Delta P_{\frac{t_{0}}{4}}^{\ell}\omega\mid\mid_{\infty}=
\mid\mid\tilde{\Delta}\tilde{P}_{\frac{t_{0}}{4}}^{\ell}\pi^{\star}\omega\mid\mid_{\infty}.
\end{equation*}
Mais d'après \cite{Ko91}, comme $\tilde{M}$ est à géométrie bornée, on a :
\beu
\mid\mid\tilde{\Delta}\tilde{P}_{\frac{t_{0}}{4}}^{\ell}(\tilde{x},\tilde{y})\mid\mid\leq ce^{-\frac{2\rho}{t_{0}}\tilde{\delta}(x,y)}
\end{equation*}
où $\tilde{\delta}$ est la distance riemannienne sur $\tilde{M}$ et $\rho$ une constante qui ne dépend que la géométrie de $\tilde{M}$. De plus,
\beu
\text{Vol}(B_{x}(r))=\mid B_{x}(r)\mid\leq ce^{\gamma r},
\end{equation*}
$\gamma$ ne d\'epend que de la g\'eom\'etrie de $\tilde{M}$. On veut voir que
\beu
\int_{\tilde{M}}e^{-\frac{\rho}{t_{0}}\tilde{\delta}(x,y)}d\sigma(x)<c\;\;\text{ et }\;\;\int_{\tilde{M}}e^{-\frac{\rho}{t_{0}}\tilde{\delta}(x,y)}d\sigma(y)<c,
\end{equation*}
où $c$ ne dépend que de $t$ et $\rho$. On note $\{x,\;e^{-\frac{\rho}{t_{0}}\tilde{\delta}(x,y)}>\alpha\}$ l'ensemble des points $x$ de $\tilde{M}$ tels que 
$e^{-\frac{\rho}{t_{0}}\tilde{\delta}(x,y)}>\alpha$. On calcule pour cela la fonction de répartition à $y$ fixé, puis à $x$ fixé, et l'on trouve :
\begin{eqnarray*}
\mid\{x,\;e^{-\frac{\rho}{t_{0}}\tilde{\delta}(x,y)}>\alpha\}\mid&\leq&c(\frac{1}{\alpha})^{\frac{\gamma}{2\rho}t_{0}},\\
\mid\{y,\;e^{-\frac{\rho}{t_{0}}\tilde{\delta}(x,y)}>\alpha\}\mid&\leq&c(\frac{1}{\alpha})^{\frac{\gamma}{2\rho}t_{0}}
\end{eqnarray*}
Si $\frac{\gamma}{2\rho}t_{0}<1$ on voit que l'on a l'inégalité souhaitée. Ce qui entraine que :
\beu
\mid\mid\tilde{\Delta}\tilde{P}_{\frac{t_{0}}{4}}^{\ell}\pi^{\star}\omega\mid\mid_{\infty}\leq c\mid\mid\omega\mid\mid_{\infty}
\end{equation*}
et
\begin{eqnarray*}
\mid\mid\Delta P_{\frac{t_{0}}{4}}^{\ell}\omega\mid\mid_{\infty}&\leq& c\mid\mid\omega\mid\mid_{\infty}\\
\mid\mid\Delta P_{\frac{t_{0}}{4}}^{\ell}\omega\mid\mid_{1}&\leq& c\mid\mid\omega\mid\mid_{1}.
\end{eqnarray*}
Pour $p=1$, on voit que :
\begin{equation}\label{eq1bis}
\mid\mid\frac{\partial}{\partial s}P^{\ell}_{s+\frac{t_{0}}{2}}\omega\mid\mid_{1}=\mid\mid\Delta P^{\ell}_{\frac{t_{0}}{2}}P^{\ell}_{s}\omega\mid\mid_{1}
\leq c\mid\mid P^{\ell}_{s}\omega\mid\mid_{1}\leq c\mid\mid \omega\mid\mid_{1}.
\end{equation}
Si $p=2$, on a vu pr\'ec\'edemment (\ref{eq1}) que 
\begin{equation}\label{eq2bis}
\mid\mid\frac{\partial}{\partial s}P^{\ell}_{s}\omega\mid\mid_{2}\leq ce^{-\tau t}\mid\mid\omega\mid\mid_{2}.
\end{equation}
On veut passer de $p=1$, $2$ \`a $1<p<2$ dans les in\'egalit\'es (\ref{eq1bis}) et (\ref{eq2bis}). Soit $1<p<2$, on veut estimer la norme :
$$\mid\mid\frac{\partial}{\partial s}P_{s+\frac{t_{0}}{4}}^{\ell}(\omega)\mid\mid_{p}.$$
Pour cela on fait appel au th\'eor\`eme de convexit\'e de Riesz-Thorin dont l'\'enonc\'e convenable dans ce contexte est le suivant.

{\it Th\'eor\`eme de Riesz-Thorin}\\
{\it On consid\`ere deux valeurs $1\leq p_{0},p_{1}\leq\infty$, $p_{0}\neq p_{1}$ ainsi qu'un op\'erateur $T$ sur l'espace des sections (modulo \'egalit\'e presque partout) du fibr\'e des $\ell$-formes. 
On suppose que $T$ applique $L^{p_{0}}(\Lambda^{\ell}M)$ dans lui-m\^eme avec une norme qui ne d\'epasse pas $M_{0}$ et d'autre part que $T$ applique $L^{p_{1}}(\Lambda^{\ell}M)$ 
dans lui-m\^eme avec une norme qui ne d\'epasse pas $M_{1}$. Alors pour tout $0<\theta<1$ et $p$ d\'efini par : $\frac{1}{p}=\frac{1-\theta}{p_{0}}+\frac{\theta}{p_{1}}$, $T$ applique 
$L^{p}(\Lambda^{\ell}M)$ dans lui-m\^eme avec une norme qui ne d\'epasse pas $M_{0}^{1-\theta}M_{1}^{\theta}$}.\\
On applique ce th\'eor\`eme avec $p_{0}=1$, $p_{1}=2$, $\theta=\frac{2}{p^{\prime}}$ et $T=\frac{\partial}{\partial s}P^{\ell}_{s+\frac{t_{0}}{4}}$. L'\'enonc\'e ci-dessus avec les in\'egalit\'es 
(\ref{eq1bis}) et (\ref{eq2bis}) montrent que
$$\mid\mid\frac{\partial}{\partial s}P^{\ell}_{s+\frac{t_{0}}{4}}\omega\mid\mid_{p}\leq\big( c^{\prime}e^{\alpha s}\big)^{1-\theta}\big(ce^{-\tau s}\big)^{\theta}
=c^{\prime\prime}e^{(\alpha(1-\theta)-\tau\theta)s}.$$
Il s'ensuit que si l'on prend $q_{0}=\frac{2(\alpha+\tau)}{\alpha+2\tau}$ on a $q_{0}^{\prime}=\frac{2(\alpha+\tau)}{\alpha}$ pour 
$q_{0}<p<q_{0}^{\prime}$ et 
$$\mid\mid\frac{\partial}{\partial s}P^{\ell}_{s+\frac{t_{0}}{4}}\omega\mid\mid_{p}\leq ce^{-\gamma(p)s}$$
avec $\gamma(p)=\alpha-\frac{2}{p^{\prime}}(\alpha+\tau)$. En effet par dualit\'e, il suffit de montrer l'in\'egalit\'e ci-dessus pour $q_{0}<p<2$ en prenant $\theta=2/p^{\prime}$, on voit que 
$$(1-\theta)\alpha-\tau\theta=(1-\frac{2}{p^{\prime}})\alpha-\frac{2\tau}{p^{\prime}}=\alpha-\frac{2}{p^{\prime}}(\alpha+\tau)<0.$$
Finalement, on a bien :
\begin{eqnarray*}
\mid\mid\omega-H\omega\mid\mid&=&\mid\mid\int_{0}^{\infty}\frac{\partial}{\partial s}P_{s}^{\ell}\omega ds\mid\mid\leq c\mid\mid\omega\mid\mid_{p}\\
\text{ et }\mid\mid H\omega\mid\mid_{p}&\leq& c\mid\mid\omega\mid\mid_{p}.
\end{eqnarray*}

{\bf Lemme 2.} {\em Pour tout $\omega\in C_{0}^{\infty}(\Lambda^{\ell}M)$, on pose $G\omega=\int_{0}^{\infty}P^{\ell}_{t}(1-H)\omega dt$. Il existe $p_{1}$ avec 
$1<p_{1}<2$ et $c>0$ tels que pour tout $p_{1}<p<p_{1}^{\prime}$ et pour tout $\omega\in C_{0}^{\infty}(\Lambda^{\ell}M)$, on ait :}
\beu
\mid\mid G\omega\mid\mid_{p}\leq c\mid\mid\omega\mid\mid_{p}.
\end{equation*}

{\bf Preuve du lemme :} Soit $\omega$ de classe $C^{\infty}$ à support compact. On pose $\omega_{0}=\omega-H\omega$, et l'on veut prouver que :
\beu
\mid\mid\int_{0}^{\infty}P_{t}^{\ell}\omega_{0}dt\mid\mid_{p}\leq c\mid\mid\omega\mid\mid_{p}\text{ pour }p_{1}<p<p_{2}
\end{equation*}
où $p_{i}$ sera déterminé dans la suite. En effet :
\beu
\int_{0}^{\infty}P_{t}^{\ell}\omega_{0}dt=\int_{0}^{\infty}P_{t}^{\ell}(1-H)\omega dt.
\end{equation*}
Pour $p=2$, puisque $\IP{\Delta\omega}{\omega}\geq c\mid\mid\omega\mid\mid_{2}$ sur ${\Cal H}_{2}^{\perp}$ :
\beu
\mid\mid P_{t}^{\ell}(1-H)\omega\mid\mid_{2}\leq ce^{-\tau t}\mid\mid\omega\mid\mid_{2}.
\end{equation*}
Soit $\epsilon$ suffisamment petit devant $\tau$. On pose $\displaystyle{\theta_{0}=\frac{\alpha+\epsilon}{\alpha+\tau}}$ alors $\theta_{0}<1$ et 
$$\frac{1}{p_{0}}=\frac{1-\theta_{0}}{q_{\epsilon}}+\frac{\theta_{0}}{2}\;\text{ o\`u }\; q_{\epsilon}=\frac{2(\alpha+\tau+\epsilon)}{\alpha+2\tau}.$$
On remarque, par le lemme pr\'ec\'edent, que $q_{\epsilon}>q_{0}$ et on constate que :
$$\alpha(1-\theta_{0})-\theta_{0}\tau<0.$$
Si $\displaystyle{\frac{1}{p}=\frac{1-\theta_{1}}{q_{\epsilon}}+\frac{\theta_{1}}{2}}$ avec $\theta_{0}\leq\theta_{1}<1$ alors $p_{0}\leq p<2$ et r\'eciproquement si $p_{0}<p<2$, 
$\displaystyle{\frac{1}{p}=\frac{1-\theta}{q_{\epsilon}}+\frac{\theta}{2}}$ avec $\theta_{0}<\theta\leq 1$ et $\alpha(1-\theta)-\theta\tau<0$. 

Soit $\displaystyle{\frac{1}{p}=\frac{1-\theta}{q_{\epsilon}}+\frac{\theta}{2}}$. Pour $q_{\epsilon}$ on a vu, Lemme 1 :
\begin{eqnarray*}
\mid\mid P_{t}^{\ell}(1-H)\omega\mid\mid_{q_{\epsilon}}&\leq&\mid\mid P_{t}^{\ell}\mid\mid_{q_{\epsilon}\rightarrow q_{\epsilon}}\mid\mid(1-H)\omega\mid\mid_{q_{\epsilon}\rightarrow q_{\epsilon}}\\
&\leq&\mid\mid P_{t}^{\ell}\mid\mid_{q_{\epsilon}\rightarrow q_{\epsilon}}\mid\mid\omega\mid\mid_{q_{\epsilon}}+\mid\mid H\omega\mid\mid_{q_{\epsilon}}\\
&\leq& c\mid\mid P_{t}^{\ell}\mid\mid_{q_{\epsilon}\rightarrow q_{\epsilon}}\mid\mid\omega\mid\mid_{q_{\epsilon}}\;\;\text{ (d'apr\`es le lemme 1)}\\
&\leq&ce^{\alpha t}\mid\mid\omega\mid\mid_{q_{\epsilon}}\;\;\text{ (d'apr\`es \cite{Ko91})}.
\end{eqnarray*}
Si $p=2$, $\displaystyle{\mid\mid P_{t}^{\ell}(1-H)\omega\mid\mid_{2}\leq ce^{-\tau t}\mid\mid\omega\mid\mid_{2}}$. En utilisant le th\'eor\`eme de Riesz - Thorin pr\'ec\'edent, on trouve que 
$$\mid\mid P_{t}^{\ell}(1-H)\omega\mid\mid_{p}\leq ce^{\alpha(1-\theta)t}e^{-\theta\tau t}\leq ce^{-\gamma(p)t}$$
o\`u $\gamma(p)=\theta\tau-\alpha(1-\theta)$. Comme $P_{t}^{\ell}(1-H)$ est auto-adjoint, on a :
$$\mid\mid P_{t}^{\ell}(1-H)\omega\mid\mid_{p}\leq ce^{-\gamma(p)t}\;\;\text{ pour }q_{\epsilon}<p<q_{\epsilon}^{\prime}.$$
Il en découle que 
\begin{equation}\label{star}
\mid\mid\int_{0}^{\infty}P_{t}^{\ell}(1-H)\omega dt\mid\mid_{p}\leq c\mid\mid\omega\mid\mid_{p}\text{ pour }q_{\epsilon}<p<q_{\epsilon}^{\prime},\;\text{ on prend }p_{1}=q_{\epsilon}
\end{equation}
et
\begin{equation}
\mid\mid\int_{0}^{\infty}P_{t}^{\ell}\omega_{0}dt\mid\mid_{p}\leq c\mid\mid\omega_{0}\mid\mid_{p}\leq c\mid\mid\omega\mid\mid_{p}
\end{equation}
ce qui termine la preuve du lemme. 

{\bf Fin de la d\'emonstration du th\'eor\`eme.} On voit sans peine que 
$$\omega=\Delta G\omega+H\omega=\Delta G\omega_{0}+H\omega=d\delta G\omega_{0}+\delta d G\omega_{0}+H\omega,\;\;\text{o\`u }\omega_{0}=(1-H)\omega.$$
On voit aussi facilement que 
\begin{equation}\label{stare}
\int_{0}^{\infty}P_{t}^{\ell}t^{-\frac{1}{2}}(1-H)\omega dt=\Delta^{-\frac{1}{2}}(1-H)\omega
\end{equation}
et
$$\mid\mid\Delta^{-\frac{1}{2}}(1-H)\omega\mid\mid_{p}\leq c\mid\mid\omega\mid\mid_{p}\;\;\;\text{pour }p_{1}<p<p_{2},\;\;\;(p_{1}\text{ et }p_{2}\text{ viennent du lemme 2}$$
d'apr\`es (\ref{star}).

On veut montrer que 
$$\mid\mid\delta d G\omega_{0}\mid\mid_{p}\leq c\mid\mid\omega\mid\mid_{p}\;\;\text{ et }\;\;\mid\mid d\delta G\omega_{0}\mid\mid_{p}\leq c\mid\mid\omega\mid\mid_{p}.$$
Comme $\tilde{M}$ est \`a g\'eom\'etrie born\'ee, pour toute constante $\gamma$ suffisamment grande, on a (voir \cite{Ba87} ou \cite{Lo06}) :
\beu
\mid\mid d\omega_{0}\mid\mid_{p}+\mid\mid \delta\omega_{0}\mid\mid_{p}\leq c\mid\mid(\Delta+\gamma)^{\frac{1}{2}}\omega\mid\mid_{p}.
\end{equation*}
Mais la fonction $f(x)=(1+x)^{\frac{1}{2}}-x^{\frac{1}{2}}$, pour $x>0$, est la transformée de Laplace d'une mesure bornée sur $\R^{+}$ \cite{Ba87},
par conséquent
\beu
\sqrt{\gamma}((1+\frac{x}{\gamma})^{\frac{1}{2}}-(\frac{x}{\gamma})^{\frac{1}{2}})=(x+\gamma)^{\frac{1}{2}}-x^{\frac{1}{2}}=
\sqrt{\gamma}\int_{0}^{\infty}e^{-\frac{x}{\gamma}t}d\mu(t)
\end{equation*}
et
\beu
(\Delta+\gamma)^{\frac{1}{2}}-\Delta^{\frac{1}{2}}=\sqrt{\gamma}\int_{0}^{\infty}e^{-\Delta\frac{t}{\gamma}}d\mu(t).
\end{equation*}
Il s'ensuit que :
\beu
\mid\mid(\Delta+\gamma)^{\frac{1}{2}}-\Delta^{\frac{1}{2}})\omega\mid\mid_{p}
\leq c\int_{0}^{\infty}\mid\mid e^{-\Delta\frac{t}{\gamma}}\mid\mid_{p\rightarrow p}\mid\mid\omega\mid\mid_{p}d\mu(t).
\end{equation*}
D'après le résultat antérieur
\beu
\int_{0}^{\infty}\mid\mid e^{-\Delta\frac{t}{\gamma}}\mid\mid_{p\rightarrow p}dt\leq \int_{0}^{\infty}e^{(\alpha(1-\theta)-\theta\tau)\frac{t}{\gamma}}dt<c
\end{equation*}
dans l'intervalle que nous considérons et il en résulte que :
\beu
\mid\mid(\Delta+\gamma)^{\frac{1}{2}}-\Delta^{\frac{1}{2}})\omega_{0}\mid\mid_{p}\leq c\mid\mid\omega\mid\mid_{p}.
\end{equation*}
Alors :
\beu
\mid\mid d\omega\mid\mid_{p}+\mid\mid \delta\omega\mid\mid_{p}
\leq c\mid\mid\Delta^{\frac{1}{2}}\omega\mid\mid_{p}+c\mid\mid\omega\mid\mid_{p}
\end{equation*}
et d'après (\ref{stare})
\beu
\mid\mid d\omega_{0}\mid\mid_{p}+\mid\mid \delta\omega_{0}\mid\mid_{p}
\leq c\mid\mid\Delta^{\frac{1}{2}}\omega_{0}\mid\mid_{p}.
\end{equation*}
Donc
\beu
\mid\mid d(\Delta^{-\frac{1}{2}}\omega_{0})\leq c\mid\mid\omega_{0}\mid\mid_{p}\text{ et }
\mid\mid \delta(\Delta^{-\frac{1}{2}}\omega_{0})\leq c\mid\mid\omega_{0}\mid\mid_{p}.
\end{equation*}
Puis
\beu
d\delta G\omega_{0}=d\delta\Delta^{-\frac{1}{2}}\Delta^{-\frac{1}{2}}\omega_{0}=\Delta^{-\frac{1}{2}}d\delta\Delta^{-\frac{1}{2}}\omega_{0}
\end{equation*}
car $d\delta$ commute avec $\Delta$. Par conséquent
\beu
\mid\mid d\delta G\omega_{0}\mid\mid_{p}\leq\mid\mid\Delta^{-\frac{1}{2}}d\mid\mid_{p\rightarrow p}
\mid\mid\delta\Delta^{-\frac{1}{2}}\omega_{0}\mid\mid_{p}\leq c\mid\mid\omega_{0}\mid\mid_{p}.
\end{equation*}
On a la même estimation pour
\beu
\mid\mid \delta dG\omega_{0}\mid\mid_{p}\leq c\mid\mid\omega_{0}\mid\mid_{p}.
\end{equation*}
De plus
\beu
\mid\mid\delta G\omega_{0}\mid\mid_{p}+\mid\mid dG\omega_{0}\mid\mid_{p}
\leq c\mid\mid\Delta^{\frac{1}{2}}G\omega_{0}\mid\mid_{p}=c\mid\mid\Delta^{-\frac{1}{2}}\omega_{0}\mid\mid_{p}
\leq c\mid\mid\omega_{0}\mid\mid_{p}\leq c\mid\mid\omega\mid\mid_{p}.
\end{equation*}
Considérons $p_{1}<p<p_{2}$. Si $\omega$ est $C^{\infty}$, il est clair que $H\omega$ est harmonique par la théorie $L^{2}$.
Si $\omega$ est dans $L^{p}$, la continuité de $H$ et la densité de $C^{\infty}_{0}$ dans $L^{p}$ entraine que $H\omega$ est harmonique et est
dans $L^{p}$. Soit $\omega$ dans $L^{p}$, d'après la décomposition précédente,
\beu
\omega=d\omega_{1}+\delta\omega_{2}+\omega_{3}
\end{equation*}
avec $\delta\omega_{1}\in L^{p}$, $d\omega_{2}\in L^{p}$, et $\omega_{3}\in L^{p}$ et harmonique. Si $d\omega=0$, $d\delta\omega_{2}=0$ et, comme
$\delta^{2}\omega_{2}=0$, $\delta\omega_{2}$ est harmonique, puisque $(1-H)\delta\omega_{2}=\delta\omega_{2}$, $H\delta\omega_{2}=0$,
$\frac{\partial}{\partial t}P_{t}^{\ell}\delta\omega_{2}=\Delta P_{t}^{\ell}\delta\omega_{2}=0$ et $\delta\omega_{2}=H\delta\omega_{2}=0$ car 
$H\omega=\lim_{t\rightarrow\infty}P_{t}^{\ell}\omega$, $P_{t}^{\ell}\delta\omega_{2}=\delta\omega_{2}$. Il s'en suit que $d\omega=0$ implique que $\omega=d\omega_{1}+\omega_{3}$, avec
$\omega_{1}\in L^{p}$, $\Delta\omega_{3}=0$. Les classes de cohomologie sont donc représentées par des formes harmoniques.


\section{Un corollaire.}

{\bf Corollaire :} {\em Soient $p_{1}$ et $p_{2}$ comme ci-dessus et $p_{1}<p<p_{2}$. On suppose que le rayon d'injectivit\'e de $M$ est minor\'e. 
Alors $\text{dim}({\Cal H}_{p})$ est finie si, et seulement si, $\text{dim}({\Cal H}_{2})$ est finie. De plus $\dim({\Cal H}_{p})=\dim({\Cal H}_{2})$.}

{\bf Preuve :} On suppose $p<p_{0}$ alors il est clair que ${\Cal H}_{p}\subset{\Cal H}_{p_{0}}$. En effet, soit $\omega\in{\Cal H}_{p}$; alors
$P_{t}^{\ell}\omega=\omega$ et
\beu
\mid\mid\omega\mid\mid_{p_{0}}\leq\mid\mid P_{t}^{\ell}\omega\mid\mid_{p_{0}}
\leq\mid\mid P_{t}^{\ell}\mid\mid_{p\rightarrow p_{0}}\mid\mid\omega\mid\mid_{p},
\end{equation*}
mais (voir \cite{CLY81}) :
\beu
\mid\mid P_{t}^{\ell}\mid\mid \leq ce^{-\tilde{\delta}^{2}(x,y)}=K(x,y).
\end{equation*}
Pour $x$ fixé, $K$ comme fonction de $y$ est dans $L^{q}$ avec $\frac{1}{p_{0}}=\frac{1}{p}+\frac{1}{q}-1$, avec une norme indépendante de $x$.
De même pour $y$ fixé, $K$ comme fonction de $x$ est dans le même $L^{q}$. Soit 
$$K(\mid\mid\omega\mid\mid)(x)=\int_{M}K(x,y)\mid\mid\omega_{y}\mid\mid d\sigma(y).$$
Alors
\beu
\mid\mid K(\mid\mid\omega\mid\mid)\mid\mid_{p_{0}}\leq c_{1}\mid\mid\omega\mid\mid_{p}
\end{equation*}
d'après \cite{FS74}, et $\omega_{0}$ est dans $L^{p_{0}}$, de plus :
\beu \dim({\Cal H}_{p})\leq \dim({\Cal H}_{p_{0}})\leq \dim({\Cal
H}_{2}).
\end{equation*}
L'espace dual de ${\Cal H}_{p}$ est $L^{q}(\Lambda^{\ell}M)/{\Cal H}_{p}^{\perp}$, avec $\frac{1}{p}+\frac{1}{q}=1$, où ${\Cal H}_{p}^{\perp}$ est
l'orthogonal de ${\Cal H}_{p}$ dans $L^{q}(\Lambda^{\ell}M)$. Il est de dimension finie si, et seulement si, ${\Cal H}_{p}$ est de dimension finie.
De plus, d'après la décomposition précédente, $L^{q}(\Lambda^{\ell}M)/{\Cal H}_{p}^{\perp}\simeq {\Cal H}_{q}$. Par conséquent
$\dim({\Cal H}_{p})=\dim({\Cal H}_{q})$. En supposant $p_{1}<p\leq 2$, on a les inégalités suivantes :
\beu
\dim({\Cal H}_{p})\leq\dim({\Cal H}_{2})\leq\dim({\Cal H}_{q})=\dim({\Cal H}_{p}).
\end{equation*}
Finalement,
\beu
\dim({\Cal H}_{p})=\dim({\Cal H}_{2}).
\end{equation*}

\end{document}